\documentclass[twoside,12pt]{article}
\usepackage[]{amsmath,amsthm,amssymb}
\usepackage[latin1]{inputenc}

\begin{document}

\title{{\Large\bf  Discrete Fourier-Jacobi transform}}

\author{Semyon  YAKUBOVICH}
\maketitle

\markboth{\rm \centerline{ Semyon   YAKUBOVICH}}{}
\markright{\rm \centerline{FOURIER-JACOBI TRANSFORM}}

\begin{abstract} {\noindent Discrete analogs of the classical Fourier-Jacobi transform are introduced and investigated. It involves series and integrals with respect to parameters of the Gauss hypergeometric function ${}_2F_1(a+in/2,a-in/2;\  c; -x^2 ), \ x >0, n \in \mathbb{N}, a,c > 0, i $ is the imaginary unit. The corresponding inversion formulas for suitable  functions and sequences in terms of these series and integrals are established.  }

\end{abstract}
\vspace{4mm}



\vspace{4mm}

\section {Introduction and preliminary results}

As is known (cf. [1]), the Fourier-Jacobi or Jacobi transform of an arbitrary function $f$ with respect to parameters of the Gauss hypergeometric function  can be defined by the following integral

$$F(\tau)=  \int_0^\infty  {}_2F_1\left(a+ {i \tau\over 2},\ a- {i\tau\over 2};\  c; -x^2 \right) f(x) dx,\  \tau >0,\eqno(1.1)$$
where $a, c >0,\ \Gamma(z)$ is the Euler gamma function [3], Vol. III.  It was introduced by Olevskii [2] in 1949 who proved the inversion formula in the space of integrable functions, having the form (cf.  [4], formula (7.51))

$$f(x)= { x^{2c-1} \over 2\pi  } \int_0^\infty  \left|{\Gamma(c-a+i\tau/2)\Gamma(a+i\tau/2) \over \Gamma(c) \Gamma(i\tau)} \right|^2  {}_2F_1\left(c-a+{i\tau\over 2},\ c-a-{i\tau\over 2};\  c; -x^2 \right)$$

$$\times  F(\tau) d\tau,\ x >0.\eqno(1.2)$$
This is why operator (1.1) is also called the Olevskii transform [4],  [5], or the index hypergeometric transform, or the  $ {}_2F_1$-index transform. The transformation (1.1) represents a class of the so-called index transforms which were  investigated by the author in [5]. 

The main goal of this paper is to investigate the  discrete analogs of the Fourier-Jacobi transform (1.1), (1.2)

$$f(x)= \sum_{n=1}^\infty  a_n\  {}_2F_1\left(a+ {in\over 2},\ a- {in\over 2};\  c; -x^2 \right),\quad x >0,\eqno(1.3)$$

$$a_n=  \int_0^\infty  {}_2F_1\left(a+ {in\over 2},\ a- {in\over 2};\  c; -x^2 \right) f(x) dx,\  n \in \mathbb{N},\eqno(1.4)$$
for suitable function $f$ and sequence $\{a_n\}_{n\ge 1}$.  Our approach suggests to use the  classical Fourier series and some integrals, involving the Gauss hypergeometric function.  
In fact, the crucial result in our investigation  will be  the following lemma. 

{\bf Lemma 1}. {\it Let  $a > 0,\  c > \hbox{max}\ (0, a- 1/2),\ u \in \mathbb{R},\ n \in \mathbb{N}$. Then the  equality holds}

$$ \Gamma(2(c-a) +1) \left|{\Gamma\left(a+ in/2\right)\over \Gamma(c)} \right|^2 \int_0^\infty y^{c-1}  {}_2F_1\left(c-a+ {1\over 2},\  c-a +1; \ c;\  - {y\over \cosh^2(u)} \right)$$

$$\times   {}_2F_1\left(a+ {in\over 2},\ a- {in\over 2};\  c; -y \right)dy =  { 2^{2(c-a)+1} \pi \sin( nu) \   [\cosh(u)]^{ 2(c-a) +1} \over \sinh (u) \sinh(\pi n)}.\eqno(1.5)$$

\begin{proof}  Since the Gauss hypergeometric functions in (1.5) behave as (see [3], Vol. III )

$$ {}_2F_1\left(a+ {in\over 2},\ a- {in\over 2};\  c; -y \right) = O(y^{-a}),\ y \to \infty,\eqno(1.6)$$

$${}_2F_1\left(c-a+ {1\over 2},\  c-a +1; \ c;\  - {y\over \cosh^2(u)} \right) = O(y^{a-c-1/2}),\ y \to \infty,\eqno(1.7)$$
integral (1.5) converges absolutely.  Then, appealing to the Mellin-Barnes representation of the function (1.7) (see [3], Vol. III, Entry 8.4.49.13)

$$ x^{c-1}  {}_2F_1\left(c-a+ {1\over 2},\  c-a +1; \ c;\  - {x\over \cosh^2(u)} \right) = {4^{c-a} \   [\cosh(u)]^{ 2(c-1)} \Gamma(c)\over \Gamma( 2(c-a)+1)\  2\pi^{3/2} i} $$

$$\times \int_{\gamma-i\infty}^{\gamma+i\infty} \Gamma(3/2-a- s)\Gamma(2-a-s) {\Gamma(s+c-1)\over \Gamma(1-s)} \left( {x\over \cosh^2(u)}\right)^{-s} ds,\eqno(1.8)$$
where  $\hbox{max} (1-c, 1-a)  < \gamma < \hbox{min} (1, 3/2- a)$,  we substitute the right-hand side of (1.8) into (1.5) and change the order of integration by Fubini's theorem due to the estimate

$$ \int_0^\infty y^{-\gamma}  \left| {}_2F_1\left(a+ {in\over 2},\ a- {in\over 2};\  c; -y \right) \right| dy$$

$$\times  \int_{\gamma-i\infty}^{\gamma+i\infty} \left| \Gamma(3/2-a- s)\Gamma(2-a-s) {\Gamma(s+c-1)\over \Gamma(1-s)} ds\right| < \infty.$$
This  estimate, indeed,  holds via (1.6) and the Stirling asymptotic formula for the gamma function [4].  Consequently, recalling formula  in [3], Vol. III, Entry 8.4.49.13 and the duplication formula for the gamma function, we derive 

$$ \Gamma(2(c-a) +1) \left|{\Gamma\left(a+ in/2\right)\over \Gamma(c)} \right|^2 \int_0^\infty y^{c-1}  {}_2F_1\left(c-a+ {1\over 2},\  c-a +1; \ c;\  - {y\over \cosh^2(u)} \right)$$

$$\times   {}_2F_1\left(a+ {in\over 2},\ a- {in\over 2};\  c; -y \right)dy =  {4^{c- 3/2}\over 2\pi i} \  [\cosh(u)]^{ 2(c-1)}    \int_{2\gamma-i\infty}^{2\gamma+i\infty} \Gamma(3- 2a- s) $$

$$\times \Gamma\left({s+ 2a - 2+ in\over 2} \right) \Gamma\left({s+ 2a-  2 - in\over 2} \right)\left(  2\cosh (u)\right)^{s} ds$$

$$= {4^{c-a- 1/2}\over 2\pi i} \  [\cosh(u)]^{ 2(c-a) } \int_{\mu-i\infty}^{\mu+i\infty} \Gamma(1- s)  \Gamma\left({s+ in\over 2} \right) \Gamma\left({s - in\over 2} \right)\left(  2\cosh (u)\right)^{s} ds,\eqno(1.9)$$
where $\mu= 2(\gamma+ a-1)$. The latter integral can be written, employing the Parseval equality for the Mellin transform [4] and formulas   in [3], Vol. III, Entries  8.4.3.1, 8.4.23.1.  Hence we get from (1.9)

$$ \Gamma(2(c-a) +1) \left|{\Gamma\left(a+ in/2\right)\over \Gamma(c)} \right|^2 \int_0^\infty y^{c-1}  {}_2F_1\left(c-a+ {1\over 2},\  c-a +1; \ c;\  - {y\over \cosh^2(u)} \right)$$

$$\times   {}_2F_1\left(a+ {in\over 2},\ a- {in\over 2};\  c; -y \right)dy $$

$$= 4^{c-a +1/2} \  [\cosh(u)]^{ 2(c-a) }\int_0^\infty e^{-x} K_{in}\left({x\over \cosh(u)}\right) dx,\eqno(1.10)$$
where $K_\nu(z)$ is the modified Bessel function [5].   But the  integral on the right-hand side of (1.10) is calculated via  [3], Vol. II, Entry 2.16.6.1

$$ \int_0^\infty  e^{-x\cosh(u)} K_{in}(x)  dx = {\pi \sin( nu) \over \sinh (u) \sinh(\pi n)}.\eqno(1.11)$$
Therefore, making a simple change of variables,  we end up with (1.5), completing the proof of Lemma 1. 

\end{proof}

\section{Inversion theorems}

We begin with

{\bf Theorem 1.} {\it Let   $a > 0,\   \hbox{max}\ (1/2, 2a- 1/2) < c < 2a +1/2$ and $0 < a < 1/2$ when $c= 1/2$. If the sequence $\{a_n\}_{n\ge 1}$ satisfies the condition 

$$ \sum_{m=1}^\infty {\left| a_m \right| e^{-\delta m} \over \left|\Gamma\left(a+ im/2\right)\right|^2 }   < \infty,\quad  \delta \in \left[0, \  {\pi\over 2} \right),\eqno(2.1)$$
then the discrete transformation $(1.3)$ can be inverted by the formula

$$a_n= {4^{a-c} \over \pi^2}  \Gamma(2(c-a) +1) \left|{\Gamma\left(a+ in/2\right)\over \Gamma(c)} \right|^2 \sinh(\pi n) \int_0^\infty x^{2c-1} \Phi_n(x) f(x) dx,\eqno(2.2)$$
where

$$  \Phi_n(x) =  \int_{-\pi}^\pi  {}_2F_1\left(c-a+ {1\over 2},\  c-a +1; \ c;\  - {x^2\over \cosh^2(u)} \right) { \tanh(u)  \sin( nu) \over [\cosh(u)]^{ 2(c-a) } } du\eqno(2.3)$$
and integral $(2.2)$ converges absolutely.}

\begin{proof}  In order to proceed the derivation of the inversion formula (2.3) we will need an  estimation of the Gauss hypergeometric function.  To do this,  we employ its integral representation in terms of the product of Bessel functions (see [3], Vol. II, Entry 2.16.21.1)

$${}_2F_1\left(a+ {in\over 2},\ a- {in\over 2};\  c; -x^2 \right) = {2^{1+c-2a} x^{1-c} \Gamma(c) \over  \left|\Gamma\left(a+ in/2\right)\right|^2 }\int_0^\infty y^{2a-c} J_{c-1} (xy) K_{in}(y) dy.\eqno(2.4)$$
 Then, appealing to the inequality  for the modified Bessel function (see [5], p. 15)

$$ \left| K_{i\tau} (x) \right| \le e^{-\delta \tau} K_0\left( x \cos(\delta) \right), \quad x, \tau >0,\ \delta \in \  \left[0, \  {\pi\over 2} \right)\eqno(2.5)$$
and the inequality $x^{1/2} |J_\nu(x)| < C_\nu,\ \nu \ge - 1/2$, where $C_\nu > 0$ is a constant,  we deduce

$$\left| {}_2F_1\left(a+ {in\over 2},\ a- {in\over 2};\  c; -x^2 \right) \right| \le {2^{1+c-2a} x^{1-c} e^{-\delta n} \Gamma(c) \over  \left|\Gamma\left(a+ in/2\right)\right|^2 }$$

$$\times \int_0^\infty y^{2a-c} \left|J_{c-1} (xy)\right|  K_{0}\left(y \cos(\delta)\right) dy $$

$$\le  {2^{1+c-2a} x^{1/2 -c} e^{-\delta n} \Gamma(c) A_c \over  \left|\Gamma\left(a+ in/2\right)\right|^2 }\int_0^\infty y^{2a-c-1/2} K_{0}\left(y \cos(\delta)\right) dy $$

$$=  {\Gamma(c) A_c\ x^{1/2 -c}  e^{-\delta n} \over  \sqrt 2\   [\cos(\delta)]^{2a-c+1/2}}  \left| {\Gamma\left( a-  c/2  + 1/4 \right)\over \Gamma\left(a+ in/2\right)}\right|^2,\quad {1\over 2} \le c < 2a + {1\over 2},\eqno(2.6)$$
where $A_c > 0$ is a constant.   Thus, substituting the value of $f$ by formula (1.3) on the right-hand side of (2.2), we change the order of integration and summation due to the assumption (2.1), inequality (2.6) and the estimate 

$$ \int_0^\infty x^{c-1/2}  \int_{0}^\pi  \left|{}_2F_1\left(c-a+ {1\over 2},\  c-a +1; \ c;\  - {x^2\over \cosh^2(u)} \right) \right| $$

$$\times { \sinh(u) \over [\cosh(u)]^{ 2(c-a) +1} }  du dx  \sum_{m=1}^\infty {\left| a_m \right| e^{-\delta m} \over \left|\Gamma\left(a+ im/2\right)\right|^2 } < \infty\eqno(2.7)$$
under the condition $c > 2a- 1/2$ (see (1.7)).  Consequently,  appealing to Lemma 1,  we obtain

$$  {4^{a-c} \over \pi^2}  \Gamma(2(c-a) +1) \left|{\Gamma\left(a+ in/2\right)\over \Gamma(c)} \right|^2 \sinh(\pi n) \int_0^\infty x^{2c-1} \Phi_n(x) f(x) dx $$

$$ =  {1\over \pi} \left|\Gamma\left(a+ in/2\right) \right|^2 \sinh(\pi n) \int_{-\pi}^\pi  \sin( nu)  \sum_{m=1}^\infty  a_m\ { \sin( mu) \over |\Gamma\left(a+ im/2\right) |^2 \sinh(\pi m)} du = a_n.$$
Theorem 1 is proved. 

\end{proof}

Concerning Fourier-Jacobi transform (1.4), we have the following result.

{\bf Theorem 2}.   {\it Let parameters $a, c$ satisfy conditions of Theorem $1$ and $f$ be a complex-valued function on $\mathbb{R}_+$ which is represented by the integral 

$$f(x) =  x^{2c-1}  \int_{-\pi}^\pi   {}_2F_1\left(c-a+ {1\over 2},\  c-a +1; \ c;\  - {x^2\over \cosh^2(u)} \right) {\varphi(u) \over [\cosh(u)]^{2(c-a)+1} } du,\quad x >0,\eqno(2.8)$$ 
where $ \varphi(u) = \psi(u)\sinh(u)$ and $\psi$ is a  $2\pi$-periodic function, satisfying the Lipschitz condition on $[-\pi, \pi]$, i.e.

$$\left| \psi(u) - \psi(v)\right| \le C |u-v|, \quad  \forall \  u, v \in  [-\pi, \pi],\eqno(2.9)$$
where $C >0$ is an absolute constant.  Then the following inversion formula for  transformation $(1.4)$  holds

$$ f(x)  =  { 4^{a-c - 1/2} \Gamma(2(c-a) +1) x^{2c-1} \over  [\pi \Gamma(c)]^2} \   \sum_{n=1}^\infty    \sinh(\pi n)   \left|\Gamma\left(a+ in/2\right)\right|^2 \Phi_n (x) a_n,\quad x > 0,\eqno(2.10)$$
where $\Phi_n (x)$ is defined by $(2.3)$.}

\begin{proof}  The proof is based on the integral for the product of the Bessel and Gauss hypergeometric functions (see [3], Vol. III, Entry 2.21.4.2)

$$\int_0^\infty y^c J_{c-1}(xy)\  {}_2F_1\left(c-a+ {1\over 2},\  c-a +1; \ c;\  - {y^2\over \cosh^2(u)} \right) dy$$

$$ = { 2^{c} x^{c-2a}   \Gamma(c) \  [\cosh(u)]^{2(c-a)+1}\  e^{-x\cosh(u)} \over  \Gamma(2(c-a) +1) },\quad x > 0,\ c > \hbox{max} \left(0,\  2a - {3\over 2}\right).\eqno(2.11)$$ 
In fact, pugging the right hand-side of (2.4) in (1.4), we change the order of integration to get

$$a_n=   {2^{1+c-2a}  \Gamma(c) \over  \left|\Gamma\left(a+ in/2\right)\right|^2 } \int_0^\infty y^{2a-c}  K_{in}(y)  \int_0^\infty x^{1-c}  J_{c-1} (xy) f(x) dx dy.\eqno(2.12)$$
The interchange follows via Fubini's theorem by virtue of the estimate (see (2.5), (2.7), (2.8) and relation 2.16.2.2 in [3], Vol. II)

$$ \int_0^\infty y^{2a-c}  \left|K_{in}(y) \right|  \int_0^\infty x^{1-c}  \left|J_{c-1} (xy) f(x) \right|dx dy $$

$$\le 2 \int_0^\infty y^{2a-c-1/2} K_{0}(y) dy   \int_0^\infty x^{c-1/2}   \int_{0}^\pi  { |\varphi(u)| \over [\cosh(u)]^{2(c-a)+1} }$$

$$\times   \left| {}_2F_1\left(c-a+ {1\over 2},\  c-a +1; \ c;\  - {x^2\over \cosh^2(u)} \right) \right| du dx$$

$$=  2^{2a-c-1/2} \Gamma^2 \left({2a-c+1/2\over 2}\right)  \int_0^\infty x^{c-1/2}  \left| {}_2F_1\left(c-a+ {1\over 2},\  c-a +1; \ c;\  - x^2 \right) \right| dx$$

$$\times  \int_{0}^\pi  { |\varphi(u)|  du \over [\cosh(u)]^{1/2+c-2a} } < \infty.$$
Hence, returning to (2.12), we substitute $f(x)$ by formula (2.8) and using  values of  integrals (1.11), (2.11), it becomes

$$ a_n= {2^{1+2(c-a)}  \pi \Gamma^2 (c) \over  \Gamma(2(c-a) +1)\  \sinh(\pi n) \left|\Gamma\left(a+ in/2\right)\right|^2 } \int_{-\pi}^\pi  \varphi(u) \  { \sin( nu) \over \sinh (u)}du.\eqno(2.13)$$
Therefore, following the same scheme as in the proof of Theorem 5 in [6],  we substitute the value of $a_n$ by (2.13) and $\Phi_n (x)$ by (2.3) into the partial sum of the series (2.10). Then,    calculating this sum via the known identity,  we invoke  the definition of $\varphi$ to obtain 

$$S_N(x)= { 4^{a-c - 1/2} \Gamma(2(c-a) +1) x^{2c-1} \over  [\pi \Gamma(c)]^2} \   \sum_{n=1}^N    \sinh(\pi n)   \left|\Gamma\left(a+ in/2\right)\right|^2 \Phi_n (x) a_n$$

$$ = { 1 \over \pi }  \sum_{n=1}^N  \int_{-\pi}^\pi  {}_2F_1\left(c-a+ {1\over 2},\  c-a +1; \ c;\  - {x^2\over \cosh^2(t)} \right) { \tanh(t)  \sin( nt) \over [\cosh(t)]^{ 2(c-a) } } dt $$

$$\times   \int_{-\pi}^\pi  \varphi(u) \  { \sin( nu)  \over \sinh (u)}du$$

$$ = { x^{2c-1}  \over 4\pi }  \int_{-\pi}^\pi  {}_2F_1\left(c-a+ {1\over 2},\  c-a +1; \ c;\  - {x^2\over \cosh^2(t)} \right) { \tanh(t)   \over [\cosh(t)]^{ 2(c-a) } }  $$

$$\times   \int_{-\pi}^\pi  \left[ \psi(u)- \psi(-u) \right]  \  {\sin \left((2N+1) (u-t)/2 \right)\over \sin( (u-t) /2)}  du dt.\eqno(2.14)$$
Since $\psi$ is $2\pi$-periodic, we treat  the latter integral with respect to $u$ as follows 

$$  \int_{-\pi}^{\pi}  \left[ \psi(u)- \psi(-u) \right]  \  {\sin \left((2N+1) (u-t)/2 \right)\over \sin( (u-t) /2)}  du $$

$$=  \int_{ t-\pi}^{t+ \pi}  \left[ \psi(u)- \psi(-u) \right]  \  {\sin \left((2N+1) (u-t)/2 \right)\over \sin( (u-t) /2)}  du $$

$$=  \int_{ -\pi}^{\pi}  \left[ \psi(u+t)- \psi(-u-t) \right]  \  {\sin \left((2N+1) u/2 \right)\over \sin( u /2)}  du. $$
Moreover,

$$ {1\over 2\pi} \int_{ -\pi}^{\pi}  \left[ \psi(u+t)- \psi(-u-t) \right]  \  {\sin \left((2N+1) u/2 \right)\over \sin( u /2)}  du - \left[ \psi(t)- \psi(-t) \right] $$

$$=  {1\over 2\pi} \int_{ -\pi}^{\pi}  \left[ \psi(u+t)- \psi(t) + \psi (-t) - \psi(-u-t) \right]  \  {\sin \left((2N+1) u/2 \right)\over \sin( u /2)}  du.$$
When  $u+t > \pi$ or  $u+t < -\pi$ then we interpret  the value  $\psi(u+t)- \psi(t)$ by  formulas

$$\psi(u+t)- \psi(t) = \psi(u+t-2\pi)- \psi(t - 2\pi),$$ 

$$\psi(u+t)- \psi(t) = \psi(u+t+ 2\pi)- \psi(t +2\pi),$$ 
respectively.  Analogously, the value  $\psi(-u-t)- \psi(-t)$  can be treated.   Then   due to the Lipschitz condition (2.9) we have the uniform estimate
for any $t \in [-\pi,\pi]$

$${\left|  \psi(u+t)- \psi(t) + \psi (-t) - \psi(-u-t) \right| \over | \sin( u /2) |}  \le 2C \left| {u\over \sin( u /2)} \right|.$$
Therefore,  owing to the Riemann-Lebesgue lemma

$$\lim_{N\to \infty } {1\over 2\pi} \int_{ -\pi}^{\pi}  \left[ \psi(u+t)- \psi(-u-t)  - \psi(t) + \psi (-t) \right]  \  {\sin \left((2N+1) u/2 \right)\over \sin( u /2)}  du =  0\eqno(2.15)$$
for all $ t\in [-\pi,\pi].$    Besides, returning to (2.14), we estimate the iterated integral 

$$ \int_{-\pi}^\pi  \left| {}_2F_1\left(c-a+ {1\over 2},\  c-a +1; \ c;\  - {x^2\over \cosh^2(t)} \right) \right| { | \tanh(t)|   \over [\cosh(t)]^{ 2(c-a) } }  $$

$$\times \int_{ -\pi}^{\pi} \left| \left[ \psi(u+t)- \psi(-u-t)  - \psi(t) + \psi (-t) \right]  \  {\sin \left((2N+1) u/2 \right)\over \sin( u /2)}  \right| du dt$$

$$\le  4 C \int_{0}^\pi  \left| {}_2F_1\left(c-a+ {1\over 2},\  c-a +1; \ c;\  - {x^2\over \cosh^2(t)} \right) \right| {  \tanh(t)  \over [\cosh(t)]^{ 2(c-a) } }  dt $$

$$\times   \int_{ -\pi}^{\pi}   \left| {u\over \sin( u /2)} \right| du < \infty,\ x >0.$$
Consequently, via  the dominated convergence theorem it is possible to pass to the limit when $N \to \infty$ under the  integral sign, and recalling (2.15), we derive

$$  \lim_{N \to \infty}   { x^{2c-1}\over 4\pi }  \int_{-\pi}^\pi  {}_2F_1\left(c-a+ {1\over 2},\  c-a +1; \ c;\  - {x^2\over \cosh^2(t)} \right) { \tanh(t)   \over [\cosh(t)]^{ 2(c-a) } }  $$

$$\times  \int_{ -\pi}^{\pi}  \left[ \psi(u+t)- \psi(-u-t)  - \psi(t) + \psi (-t) \right] $$

$$\times  \  {\sin \left((2N+1) u/2 \right)\over \sin( u /2)}  du dt $$

$$=  {x^{2c-1}  \over 4\pi }  \int_{-\pi}^\pi  {}_2F_1\left(c-a+ {1\over 2},\  c-a +1; \ c;\  - {x^2\over \cosh^2(t)} \right) { \tanh(t)   \over [\cosh(t)]^{ 2(c-a) } }    $$

$$ \times \lim_{N \to \infty}  \int_{ -\pi}^{\pi}  \left[ \psi(u+t)- \psi(-u-t)  - \psi(t) + \psi (-t) \right]  \  {\sin \left((2N+1) u/2 \right)\over \sin( u /2)}  du dt = 0.$$
Hence, combining with (2.14),  we obtain  by virtue of  the definition of $\varphi$ and $f$

$$ \lim_{N \to \infty}  S_N(x) =   {x^{2c-1} \over 2}    \int_{-\pi}^\pi   {}_2F_1\left(c-a+ {1\over 2},\  c-a +1; \ c;\  - {x^2\over \cosh^2(t)} \right)$$

$$\times  {  \varphi (t)+ \varphi(-t)  \over [\cosh(t)]^{ 2(c-a)+1 } }dt = f(x),$$
where the integral (2.8) converges since $\varphi \in C[-\pi,\pi]$.  Thus we established  (2.10), completing the proof of Theorem 2.
 
\end{proof}

\bigskip
\centerline{{\bf Acknowledgments}}
\bigskip

\noindent The work was partially supported by CMUP, which is financed by national funds through FCT (Portugal)  under the project with reference UIDB/00144/2020.

\bigskip
\centerline{{\bf References}}
\bigskip
\baselineskip=12pt
\medskip
\begin{enumerate}

\item[{\bf 1.}\ ]  T. Koornwinder, A new proof of a Paley-Wiener type theorem for the Jacobi transform,   {\it Arkiv for Matematik} {\bf 13} (1975),  145-159.

\item[{\bf 2.}\ ]  M.N. Olevskii, On the representation of an arbitrary function by integral with the kernel involving the hypergeometric function,   {\it Dokl. AN SSSR} {\bf 69} (1949), N 1, 11-14 (in Russian). 

\item[{\bf 3.}\ ] A.P. Prudnikov, Yu.A. Brychkov and O.I. Marichev, {\it Integrals and Series}. Vol. I: {\it Elementary
Functions}, Vol. II: {\it Special Functions}, Gordon and Breach, New York and London, 1986, Vol. III : {\it More special functions},  Gordon and Breach, New York and London,  1990.

\item[{\bf 4.}\ ]  S. Yakubovich and Yu. Luchko, The Hypergeometric Approach to Integral Transforms and Convolutions, {\it Kluwer
Academic Publishers, Mathematics and Applications.} Vol.287, 1994.

\item[{\bf 5.}\ ] S. Yakubovich, {\it Index Transforms}, World Scientific Publishing Company, Singapore, New Jersey, London and
Hong Kong, 1996.

\item[{\bf 6.}\ ]  S. Yakubovich, Discrete Kontorovich-Lebedev transforms, {\it The Ramanujan Journal} (to appear).

\end{enumerate}

\vspace{5mm}

\noindent S.Yakubovich\\
Department of  Mathematics,\\
Faculty of Sciences,\\
University of Porto,\\
Campo Alegre st., 687\\
4169-007 Porto\\
Portugal\\
E-Mail: syakubov@fc.up.pt\\

\end{document}